\documentstyle[11pt, amsfonts]{amsart}
\numberwithin{equation}{section}

\newtheorem{theorem}{Theorem}

\newtheorem{prop}{Proposition}

\newtheorem{corollary}{Corollary}

\newtheorem{lemma}{Lemma}

\newcommand{\beginproof}{\noindent{\bf Proof. }}

\newcommand{\SP}{{\mathbb S}}
\newcommand{\RRR}{{\mathcal R}}

\def\R{{\mathbb R}}
\def\N{{\mathbb N}}
\def\C{{\mathbb C}}

\def\vol{\mbox{\rm vol}}

\def\endproof{\begin{flushright}
$ \Box $ \\
\end{flushright}}

\begin{document}
\hfill\today
\bigskip
\author{Alexander Koldobsky and Artem Zvavitch}\thanks{The first named author is supported in part by U.S. National Science Foundation Grant DMS-1265155. The second  named author is supported in part by U.S. National Science
Foundation Grant DMS-1101636}

\address{Department of Mathematics, University of Missouri, Columbia, MO 65211, USA} \email{koldobskiya@@missouri.edu}

\address{Department of Mathematics, Kent State University,
Kent, OH 44242, USA} \email{zvavitch@@math.kent.edu}

\title{An isomorphic version of the Busemann-Petty problem for arbitrary measures.}

\begin{abstract}
The Busemann-Petty problem for an arbitrary measure $\mu$ with non-negative even continuous density in $\R^n$ 
asks whether 
origin-symmetric convex bodies in $\R^n$ with smaller $(n-1)$-dimensional measure $\mu$ of all central
hyperplane sections necessarily have smaller measure $\mu.$ It was shown in \cite{Z2} that the answer to this problem is 
affirmative for $n\le 4$ and negative for $n\ge 5$. In this paper we prove an isomorphic version of this result.
Namely, if $K,M$ are origin-symmetric convex bodies in $\R^n$ such that 
$\mu(K\cap \xi^\bot)\le \mu(M\cap \xi^\bot)$ for every $\xi \in \SP^{n-1},$ then $\mu(K)\le \sqrt{n}\ \mu(M).$
Here $\xi^\bot$ is the central hyperplane perpendicular to $\xi.$   We also study  the above question with additional assumptions on the body $K$ and present the complex version of the problem.
In the special case where the measure $\mu$ is convex  we show that $\sqrt{n}$ can be replaced by $cL_n,$
where $L_n$ is the maximal isotropic constant. Note that, by a recent result of Klartag, $L_n \le O(n^{1/4})$. 
Finally we prove a slicing inequality
$$
\mu(K)\le C n^{1/4} \max_{\xi \in \SP^{n-1}} \mu(K \cap \xi^\perp)\ \vol_n(K)^{\frac 1n}
$$ 
for any convex even measure $\mu$ and any symmetric convex body $K$ in $\R^n,$ where $C$
is an absolute constant. This inequality was recently proved in \cite{K2} for arbitrary measures with
continuous density, but
with $\sqrt{n}$ in place of $n^{1/4}.$

\end{abstract}

\maketitle

\section{Introduction}

Let $f$ be a non-negative even continuous function on $\R^n$, and let $\mu$ be the measure in $\R^n$ with  density
$f$, i.e. for every compact set $B\subset \R^n$
$$
\mu(B)=\int\limits_B f(x)dx.
$$
This definition also applies to compact sets $B\subset \xi^\bot,$ where
$\xi\in \SP^{n-1}$ and $\xi^\bot$ is the central hyperplane
orthogonal to $\xi$. The following problem was solved in \cite{Z2}.
\smallbreak

\noindent{\bf Busemann-Petty problem for general measures
(BPGM): } \nopagebreak \noindent {\it Fix $n\ge 2$. Given two
convex origin-symmetric bodies $K$ and $M$ in $\R^n$ such
that
$$
\mu(K\cap \xi^\bot)\le\mu(M\cap \xi^\bot)
$$
for every $\xi \in \SP^{n-1}$, does it follow that
$$
\mu(K)\le \mu(M)?
$$}
The BPGM problem is a triviality for $n=2$ and strictly positive $f$,
 and the answer is ``yes'',  moreover $K \subseteq M$. It was proved in \cite{Z2}, that 
 for every strictly positive density $f$ the answer to BPGM is
 affirmative if $n \le 4$ and negative if $n \ge 5$.
 
The BPGM problem is a generalization of the original Busemann-Petty  problem,
posed in 1956 (see \cite{BP}) and   asking the same question for
 Lebesgue measure $\mu(K)=\vol_n(K)$; see \cite{Zh, GKS, Ga, Kbook}  
 for the solution and historical details.
 
Since the answer to BPGM is negative in most dimensions, it is natural to consider the following question.

 \noindent{\bf Isomorphic Busemann-Petty problem for general measures: } \nopagebreak \noindent {\it Does there exist a universal constant  ${\cal L}$ such that for  any measure $\mu$ with continuous non-negative even density $f$ and any  two
origin-symmetric convex bodies $K$ and $M$ in $\R^n$ such
that
$$
\mu(K\cap \xi^\bot)\le\mu(M\cap \xi^\bot)
$$
for every $\xi \in \SP^{n-1}$, one necessarily has
$$
\mu(K)\le  {\cal L} \mu(M)?
$$}
In Section \ref{isomorphic} we give an answer to this question with a constant not depending
on the measure or bodies, but dependent on the dimension, namely  we show that one can take ${\cal L}=\sqrt{n}.$ 
We do not know whether this constant is optimal for general measures, but we are able to
improve the constant $\sqrt{n}$ to $Cn^{1/4}$ for convex measures using the techniques of Ball \cite{Ba2} and Bobkov \cite{Bob}; see Section \ref{logconcave}. We also (see the end of  Section \ref{isomorphic})  provide
better estimates under additional assumptions that $K$ is a convex $k$-intersection body or $K$ is the unit ball
of a subspace of $L_p.$ Finally,  Section \ref{complex}  is dedicated to the complex version of the isomorphic 
Busemann-Petty problem for arbitrary measures.

In the case of volume the isomorphic Busemann-Petty problem is closely related to the hyperplane problem
of Bourgain \cite{Bo1, Bo2, Bo3} which asks whether there exists 
an absolute constant $C$ so that for any origin-symmetric convex body $K$ in $\R^n$
$$
\vol_n(K)^{\frac {n-1}n} \le C \max_{\xi \in \SP^{n-1}} \vol_{n-1}(K\cap \xi^\bot);
$$
see \cite{MP} or \cite{BGVV} for the connection between these two problems.
The hyperplane problem is still open, with the best-to-date estimate $C= O(n^{1/4})$ established
by Klartag \cite{Kl}, who slightly improved the previous estimate of Bourgain \cite{Bo3}. In Section \ref{logconcave},
following recent results of  Bobkov \cite{Bob}, we show that Klartag's result can be extended to all convex  measures in the following form.
There exists an absolute constant $C$ so that for every even convex measure $\mu$ on
$\R^n$ and every origin-symmetric convex body $K$ in $\R^n$
$$\mu(K)\le Cn^{1/4} \max_{\xi\in \SP^{n-1}} \mu(K\cap \xi^\bot) \vol_n(K)^{\frac 1n}.$$
Note that this inequality was proved in \cite{K2} for arbitrary measures $\mu$ with even continuous density,
but with the constant $\sqrt{n}$ in place of $n^{1/4}:$ 
\begin{equation}\label{sqrtn}
\mu(K)\le \sqrt{n} \frac{n}{n-1} c_n \max_{\xi\in \SP^{n-1}} \mu(K\cap \xi^\bot) \vol_n(K)^{\frac 1n},
\end{equation}
where $c_n= \vol_n(B_2^n)^{\frac{n-1}n}/\vol_{n-1}(B_2^{n-1}) < 1$ and $B_2^n$ is the unit Euclidean
ball in $\R^n.$ Also, for some special classes of bodies, including
unconditional bodies, $k$-intersection bodies, duals of bodies with bounded volume ratio, inequality (\ref{sqrtn})
has been proved with an absolute constant in place of $\sqrt{n}$ (see \cite{K1,K4,K9}). Versions
of (\ref{sqrtn}) for lower dimensional sections can be found in \cite{K5}.

\section{Isomorphic Busemann-Petty problem with ${\cal L}=\sqrt{n}$} \label{isomorphic}

We need several definitions and facts.
A closed bounded set $K$ in $\R^n$ is called a {\it star body}  if 
every straight line passing through the origin crosses the boundary of $K$ 
at exactly two points different from the origin, the origin is an interior point of $K,$
and the {\it Minkowski functional} 
of $K$ defined by 
$$\|x\|_K = \min\{a\ge 0:\ x\in aK\}$$
is a continuous function on $\R^n.$ 

The {\it radial function} of a star body $K$ is defined by
$$\rho_K(x) = \|x\|_K^{-1}, \qquad x\in \R^n \setminus \{0\}.$$
If $x\in \SP^{n-1}$ then $\rho_K(x)$ is the radius of $K$ in the
direction of $x.$

If $\mu$ is a measure on $K$ with even continuous density $f$, then 
\begin{equation} \label{polar-measure}
\mu(K) = \int_K f(x)\ dx = \int\limits_{\SP^{n-1}}\left(\int\limits_0^{\|\theta\|^{-1}_K} r^{n-1} f(r\theta)\ dr\right) d\theta.
\end{equation}
Putting $f=1$, one gets
\begin{equation} \label{polar-volume}
\vol_n(K)
=\frac{1}{n} \int_{\SP^{n-1}} \rho_K^n(\theta) d\theta=
\frac{1}{n} \int_{\SP^{n-1}} \|\theta\|_K^{-n} d\theta.
\end{equation}

The {\it spherical Radon transform} 
$\RRR:C(\SP^{n-1})\mapsto C(\SP^{n-1})$  
is a linear operator defined by
$$\RRR f(\xi)=\int_{\SP^{n-1}\cap \xi^\bot} f(x)\ dx,\quad \xi\in \SP^{n-1}$$
for every function $f\in C(\SP^{n-1}).$

The polar formulas (\ref{polar-measure}) and  (\ref{polar-volume}), applied to a hyperplane section of $K$, express 
volume of such a section in terms of the spherical Radon transform:

$$\mu(K\cap \xi^\bot) = \int_{K\cap \xi^\bot} f =  
\int_{\SP^{n-1}\cap \xi^\bot} \left(\int_0^{\|\theta\|_K^{-1}} r^{n-2}f(r\theta)\ dr \right)d\theta$$
\begin{equation} \label{measure=spherradon}
=\RRR\left(\int_0^{\|\cdot\|_K^{-1}} r^{n-2}f(r\ \cdot)\ dr \right)(\xi).
\end{equation}
and
\begin{equation} \label{volume=spherradon}
\vol_{n-1}(K\cap \xi^\bot) = \frac{1}{n-1} \int_{\SP^{n-1}\cap \xi^\bot} \|\theta\|_K^{-n+1}d\theta =
\frac{1}{n-1} \RRR(\|\cdot\|_K^{-n+1})(\xi).
\end{equation}

The spherical Radon 
transform is self-dual (see \cite[Lemma 1.3.3]{Gr}), namely,
for any functions $f,g\in C(\SP^{n-1})$
\begin{equation} \label{selfdual}
\int_{\SP^{n-1}} \RRR f(\xi)\ g(\xi)\ d\xi = \int_{\SP^{n-1}} f(\xi)\ \RRR g(\xi)\ d\xi.
\end{equation}
Using self-duality, one can extend the spherical Radon transform to measures. 
Let $\nu$ be a finite Borel measure on $\SP^{n-1}.$
We define the spherical Radon transform of $\nu$ as a functional $\RRR\nu$ on
the space $C(\SP^{n-1})$ acting by
$$(\RRR\nu,f)= (\nu, \RRR f)= \int_{\SP^{n-1}} \RRR f(x) d\nu(x).$$
By Riesz's characterization of continuous linear functionals on the
space $C(\SP^{n-1})$, 
$\RRR\nu$ is also a finite Borel measure on $\SP^{n-1}.$ If $\nu$ has 
continuous density $g,$ then by (\ref{selfdual}) the 
Radon transform of $\nu$ has density $\RRR g.$

The class of intersection bodies was introduced by Lutwak \cite{Lu}.
Let $K, L$ be origin-symmetric star bodies in $\R^n.$ We say that $K$ is the 
intersection body of $L$ if the radius of $K$ in every direction is 
equal to the $(n-1)$-dimensional volume of the section of $L$ by the central
hyperplane orthogonal to this direction, i.e. for every $\xi\in \SP^{n-1},$
\begin{equation} \label{intbodyofstar}
\rho_K(\xi)= \|\xi\|_K^{-1} = |L\cap \xi^\bot|.
\end{equation} 
All bodies $K$ that appear as intersection bodies of different star bodies
form {\it the class of intersection bodies of star bodies}. 

Note that the right-hand
side of (\ref{intbodyofstar}) can be written in terms of the spherical Radon transform using (\ref{volume=spherradon}):
$$\|\xi\|_K^{-1}= \frac{1}{n-1} \int_{\SP^{n-1}\cap \xi^\bot} \|\theta\|_L^{-n+1} d\theta=
\frac{1}{n-1} \RRR(\|\cdot\|_L^{-n+1})(\xi).$$
It means that a star body $K$ is 
the intersection body of a star body if and only if the function $\|\cdot\|_K^{-1}$
is the spherical Radon transform of a continuous positive function on $\SP^{n-1}.$
This allows to introduce a more general class of bodies. A star body
$L$ in $\R^n$ is called an {\it intersection body} 
if there exists a finite Borel measure \index{intersection body}
$\nu$ on the sphere $\SP^{n-1}$ so that $\|\cdot\|_L^{-1}= \RRR\nu$ as functionals on 
$C(\SP^{n-1}),$ i.e. for every continuous function $f$ on $\SP^{n-1},$
\begin{equation} \label{defintbody}
\int_{\SP^{n-1}} \|x\|_L^{-1} f(x)\ dx = \int_{\SP^{n-1}} \RRR f(x)\ d\nu(x).
\end{equation}

Intersection bodies played an essential role in the solution of the Busemann-Petty problem;  we refer the reader to 
\cite{Ga, Kbook, KoY} for more information about intersection bodies. It was proved in \cite{Z2} (Theorems 3, 4), that if $K$ is an intersection body then the answer to the BPGM is  affirmative for $K$ and any convex symmetric body $M$, whose central sections have greater $\mu$-measure than the corresponding sections of $K$.
\smallbreak

We need the following simple fact; cf. \cite[Lemma 1]{Z2}.

\begin{lemma}\label{lemma:el}
For any $\omega, a, b >0$ and any measurable function $\alpha:\R^+ \to \R^+$ we have
\begin{equation}\label{elem}
\frac{\omega}{a} \int\limits_0^{a} t^{n-1} \alpha(t)\ dt - \omega  \int\limits_{0}^{a}t^{n-2}\alpha(t)\ dt
 \le\frac{\omega}{a} \int\limits_0^{b} t^{n-1} \alpha(t)\ dt -  \omega \int\limits_{0}^{b}t^{n-2}\alpha(t)\ dt,
\end{equation}
provided all the integrals exist.
\end{lemma}
\beginproof
The desired inequality is equivalent to
$$a\int_a^b t^{n-2} \alpha(t)\ dt \le \int_a^b t^{n-1} \alpha(t)\ dt.$$
\endproof

 Denote by
$$d_{BM}(K,L)=\inf \{d>0: \exists T \in GL(n): K \subset TL \subset dK \}$$ the  Banach-Mazur distance between 
two origin-symmetric convex bodies $L$ and $K$ in $\R^n$ (see \cite[Section 3]{MS}), and let
$$
d_{I}(K)=\min \{d_{BM}(K, L): \mbox{ $L$ is an intersection body in $\R^n$}\}.
$$

\begin{theorem}\label{th:functg}  For  any measure $\mu$ with continuous, non-negative even density $f$ on $\R^n$ and any  two
convex origin-symmetric convex bodies $K, M \subset \R^n$ such
that
\begin{equation}\label{e:measure}
\mu(K\cap \xi^\bot)\le\mu(M\cap \xi^\bot),\ \qquad \forall \xi\in \SP^{n-1}
\end{equation}
we have
$$
\mu(K)\le   d_{I}(K) \mu(M).
$$
\end{theorem}
\beginproof 

 First, we use the polar formula (\ref{measure=spherradon}) to write
the condition (\ref{e:measure}) in terms of the spherical Radon transform:
\begin{equation}\label{eq1}
\RRR\left(\int_0^{\|\cdot\|_K^{-1}} r^{n-2}f(r\ \cdot)\ dr \right)(\xi) \le
\RRR\left(\int_0^{\|\cdot\|_M^{-1}} r^{n-2}f(r\ \cdot)\ dr \right)(\xi).
\end{equation}
Next we consider an intersection body $L$ such that $L\subset K \subset  d_I(K) L$ (note that a linear image of
an intersection body is again an intersection body; see for example \cite[Theorem 8.1.6]{Ga}) and integrate (\ref{eq1}) over $\SP^{n-1}$ with respect to the
measure $\nu$ corresponding to the intersection body $L$. Using (\ref{defintbody}) we get
\begin{equation}\label{eq:cool}
\int\limits_{\SP^{n-1}}\|x\|_L^{-1}\int\limits_{0}^{\|x\|^{-1}_K} t^{n-2}
 f(tx)dt\ dx 
\le \!\!\int\limits_{\SP^{n-1}}\|x\|_L^{-1} \int\limits_{0}^{\|x\|^{-1}_M} t^{n-2} f (tx)
\label{1}
 dt\ dx.
\end{equation}
Now, we apply (\ref{elem}) with $\omega =\|x\|_L^{-1}, a =\|x\|_K^{-1}, b=\|x\|_M^{-1}$ and $\alpha(t)=f(tx)$ to get
\begin{align}\label{eq:element}
\frac{\|x\|_L^{-1}}{\|x\|_K^{-1}} \int\limits_0^{\|x\|_K^{-1}} t^{n-1} & f(tx)dt -\|x\|_L^{-1}  \int\limits_{0}^{\|x\|_K^{-1}}t^{n-2} f(tx)dt \nonumber \\
 \le&\frac{\|x\|_L^{-1}}{\|x\|_K^{-1}} \int\limits_0^{\|x\|_M^{-1}} t^{n-1} f(tx)dt -\|x\|_L^{-1}  \int\limits_{0}^{\|x\|_M^{-1}}t^{n-2} f(tx)dt.
\end{align}
Integrating (\ref{eq:element}) over $\SP^{n-1}$, adding it to (\ref{eq:cool}) and using $L\subset K \subset  d_I(K) L$ we get
\begin{equation}
\int_{\SP^{n-1}}\frac{\|x\|_L^{-1}}{\|x\|_K^{-1}} \int\limits_0^{\|x\|_K^{-1}} t^{n-1} f(tx)dt\ dx 
 \le \int_{\SP^{n-1}}\frac{\|x\|_L^{-1}}{\|x\|_K^{-1}} \int\limits_0^{\|x\|_M^{-1}} t^{n-1} f(tx)dt\ dx 
\end{equation}
and
$$
\frac{1}{d_I(K)}\int_{\SP^{n-1}} \int\limits_0^{\|x\|_K^{-1}} t^{n-1} f(tx)dt\ dx 
 \le \int_{\SP^{n-1}}\int\limits_0^{\|x\|_M^{-1}} t^{n-1} f(tx)dt\ dx.
$$
The result follows from (\ref{polar-measure}).
\endproof
 It is easy to see that the Euclidean ball $B_2^n$ is an intersection body. By John's theorem (see, for example, \cite[Section 3]{MS} or \cite[Theorem 4.2.12]{Ga}), $d_{BM} (K, B_2^n) \le \sqrt{n}$ for all convex origin-symmetric bodies $K \subset \R^n$. This immediately shows that $d_I(K) \le \sqrt{n}$ for all convex origin-symmetric bodies $K \subset \R^n$. This fact together with Theorem \ref{th:functg} implies

 \begin{corollary} \label{isom} For  any measure $\mu$ with continuous non-negative even density on $\R^n$ and any  two
convex origin-symmetric convex bodies $K, M \subset \R^n$ such
that
\begin{equation}
\mu(K\cap \xi^\bot)\le\mu(M\cap \xi^\bot),\ \qquad \forall \xi\in \SP^{n-1}
\end{equation}
we have
$$
\mu(K)\le   \sqrt{n}\ \mu(M).
$$
\end{corollary} 
 
If the body $K$ in Theorem \ref{th:functg} is an intersection body, the constant ${\cal L}=1$ (see \cite{Z2}, Theorem 1);  this is an analog of the well-known
Lutwak's connection between intersection bodies and the Busemann-Petty problem. There are other special
classes of bodies for which the constant ${\cal L}$ does not depend on the dimension.

The classes of $k$-intersection bodies were introduced 
in \cite{K6,K7}. For an integer $k,\ 1\le k <n$ and star bodies $D,L$ in $\R^n,$
we say that $D$ is the $k$-intersection body of $L$ if for every $(n-k)$-dimensional
subspace $H$ of $\R^n,$ 
$$|D\cap H^\bot|= |L\cap H|,$$
where $H^\bot$ is the $k$-dimensional subspace orthogonal to $H.$
Taking the closure in the radial metric of the class of all $D$'s that appear
as $k$-intersection bodies of star bodies, we define the class of 
{\it $k$-intersection bodies} (the original definition in \cite{K6,K7} was
different; the equivalence of definitions was proved by Milman \cite{Mi}).
These classes of bodies are important for generalizations of the Busemann-Petty
problem; see \cite{Kbook}.

To estimate the Banach-Mazur distance from $k$-intersection bodies to intersection bodies,
we use two results. The first was proved in \cite{K8}; see also \cite[Theorem 4.11]{Kbook}.

\begin{prop}\label{lplq} The unit ball of any finite dimensional subspace of $L_q$ with $0<q\le 2$
is an intersection body.
\end{prop}

We also use a result from \cite{KK}; see also \cite[Theorem 6.18]{Kbook}.

\begin{prop} \label{factor} For every $k\in \N$ and every $0<q<1,$ there exists a constant
$c(k,q)$ such that for every $n\in \N,\ n>k$ and every origin-symmetric convex 
$k$-intersection body $D$ in $\R^n$ there exists an $n$-dimensional subspace
of $L_q([0,1])$ whose unit ball $L$ satisfies $L\subset D \subset c(k,q)L.$
\end{prop}

 \begin{corollary}\label{kint} Let $k\in \N.$ There exists a constant $C(k)$ such that for any $n>k,$ any measure $\mu$ with continuous non-negative even density on $\R^n,$ any convex $k$-intersection body $K$ in $\R^n$
and any origin-symmetric convex body $M \subset \R^n$, the inequalities
$$
\mu(K\cap \xi^\bot)\le\mu(M\cap \xi^\bot),\ \qquad \forall \xi\in  \SP^{n-1}
$$
imply
$$
\mu(K)\le  C(k)\ \mu(M).
$$
\end{corollary} 
\beginproof

Let $q=1/2.$ Propositions \ref{lplq} and \ref{factor} imply that $d_I(K)\le c(k,1/2)=: C(k).$
The result follows from Theorem \ref{th:functg}.
\endproof
The constant $\sqrt{n}$ in Corollary \ref{isom} can also be improved if $K$ is the unit ball of a subspace
of $L_p,\ p>2.$ For such $K$, by  a result of Lewis \cite{Le} (see also \cite{SZ} for a different proof),
we have $d_{BM}(K,B_2^n)\le n^{1/2-1/p}.$ Since $B_2^n$ is an intersection body, Theorem \ref{th:functg}
implies the following.

\begin{corollary}\label{p>2} Let $p>2,$ let $K$ be the unit ball of an $n$-dimensional subspace
of $L_p,$ and let $\mu$ be a measure with even continuous density on $\R^n.$
Suppose that $M$ is an origin-symmetric convex body in $\R^n$ so that
$$
\mu(K\cap \xi^\bot)\le\mu(M\cap \xi^\bot),\ \qquad \forall \xi\in  \SP^{n-1}.
$$
Then $$\mu(K)\le n^{1/2-1/p} \mu(M).$$
\end{corollary}
\noindent {\bf Remark.} The statements of Theorem \ref{th:functg} and Corollaries \ref{isom}, \ref{kint}, \ref{p>2} hold true
if $M$ is any star body.

\section{The complex case} \label{complex}

Origin symmetric convex bodies in $\C^n$ are the unit balls of norms on $\C^n.$
We denote by $\|\cdot\|_K$
the norm corresponding to the body $K:$
$$K=\{z\in \C^n:\ \|z\|_K\le 1\}.$$
In order to define volume, we identify $\C^n$ with $\R^{2n}$ using the standard mapping
$$\xi = (\xi_1,...,\xi_n)=(\xi_{11}+i\xi_{12},...,\xi_{n1}+i\xi_{n2})
 \mapsto  (\xi_{11},\xi_{12},...,\xi_{n1},\xi_{n2}).$$
 Since norms on $\C^n$ satisfy the equality
$$\|\lambda z\| = |\lambda|\|z\|,\quad \forall z\in \C^n,\  \forall\lambda \in \C,$$
origin symmetric complex convex bodies correspond to those origin symmetric convex bodies
$K$  in $\R^{2n}$ that are invariant
 with respect to any coordinate-wise two-dimensional rotation, namely for each $\theta\in [0,2\pi]$
 and each $\xi= (\xi_{11},\xi_{12},...,\xi_{n1},\xi_{n2})\in \R^{2n}$
  \begin{equation} \label{rotation}
  \|\xi\|_K =
 \|R_\theta(\xi_{11},\xi_{12}),...,R_\theta(\xi_{n1},\xi_{n2})\|_K,
 \end{equation}
 where $R_\theta$ stands for  the counterclockwise rotation of $\R^2$ by the angle
 $\theta$ with respect to the origin. We shall say that $K$ is a {\it complex convex body
 in $\R^{2n}$} if $K$ is a convex body and satisfies equations (\ref{rotation}). Similarly,
 complex star bodies are $R_\theta$-invariant star bodies in $\R^{2n}.$
 \medbreak
 For $\xi\in \C^n,
|\xi|=1,$ denote by
$$H_\xi = \{ z\in \C^n:\ (z,\xi)=\sum_{k=1}^n z_k\overline{\xi_k} =0\}$$
the complex hyperplane through the origin, perpendicular to $\xi.$
 Under the standard mapping from $\C^n$ to $\R^{2n}$ the hyperplane $H_\xi$ 
 turns into a $(2n-2)$-dimensional subspace of $\R^{2n}.$ 
   \medbreak
Denote by $C_c(\SP^{2n-1})$ the space of $R_\theta$-invariant continuous functions, i.e.
  continuous real-valued functions $f$ on the unit sphere $\SP^{2n-1}$ in $\R^{2n}$ satisfying 
  $f(\xi)=f(R_\theta(\xi))$ for all $\xi\in \SP^{2n-1}$ and all $\theta\in [0,2\pi].$ The {\it complex spherical
  Radon transform} is an operator ${\cal{R}}_c: C_c(\SP^{2n-1})\to C_c(\SP^{2n-1})$ defined by
  $${\cal{R}}_cf(\xi) = \int_{\SP^{2n-1}\cap H_\xi} f(x) dx.$$
 
  We say that a finite Borel measure $\nu$ on
$\SP^{2n-1}$ is $R_\theta$-invariant if for any continuous function $f$ on $\SP^{2n-1}$ and any $\theta\in [0,2\pi]$,
$$\int_{\SP^{2n-1}} f(x) d\nu(x) = \int_{\SP^{2n-1}} f(R_\theta x) d\nu(x).$$
The complex spherical Radon transform of an $R_\theta$-invariant measure $\nu$ is defined
as a functional ${\cal{R}}_c\nu$ on the space $C_c(\SP^{2n-1})$ acting by 
$$ \left({\cal{R}}_c\nu, f \right) = \int_{\SP^{2n-1}} {\cal{R}}_cf(x) d\nu(x).$$

Complex intersection bodies were introduced and studied in \cite{KPZ}. 
An origin symmetric complex star body $K$ in $\R^{2n}$ is called a {\it complex intersection body} if there
exists a finite Borel $R_\theta$-invariant measure $\nu$ on $\SP^{2n-1}$ so that
$\|\cdot\|_K^{-2}$ and ${\cal{R}}_c\nu$ are equal as functionals on $C_c(\SP^{2n-1}),$ i.e.
for any $f\in C_c(\SP^{2n-1})$
\begin{equation}\label{defcompint}
\int_{\SP^{2n-1}} \|x\|_K^{-2} f (x)\ dx = \int_{\SP^{2n-1}} {\cal{R}}_c f(\theta) d\nu(\theta).
\end{equation}
It was proved in \cite{KPZ} that an origin-symmetric complex star body $K$ in $\R^{2n}$ is
a complex intersection body if and only if the function $\|\cdot\|_K^{-2}$ represents a positive 
definite distribution on $\R^{2n}.$
\smallbreak
We need a polar formula for the measure of a complex star body $K$ in $\R^{2n}:$
\begin{equation}\label{meas-polar}
\mu(K) = \int_K f(x)\ dx = \int_{\SP^{2n-1}} \left(\int_0^{\|\theta\|_K^{-1}} r^{2n-1}f(r\theta)\ dr \right) d\theta.
\end{equation}

For every $\xi\in \SP^{2n-1},$
$$ \mu(K\cap H_\xi) = \int_{K\cap H_\xi} f(x) dx$$
$$= \int_{\SP^{2n-1}\cap H_\xi} \left(\int_0^{\|\theta\|_K^{-1}}r^{2n-3} f(r\theta) dr \right)d\theta$$
\begin{equation} \label{meas-sect-polar}
= {\cal{R}}_c \left(\int_0^{\|\cdot\|_K^{-1}} r^{2n-3} f(r\cdot)\ dr \right)(\xi),
\end{equation}

We use an elementary inequality, which is a modification of Lemma \ref{lemma:el}.

\begin{lemma} \label{elemcomp}
For any $\omega, a, b >0$ and measurable function $\alpha:\R^+ \to \R^+$ we have
$$
\frac{\omega^2}{a^2} \int\limits_0^{a} t^{2n-1} \alpha(t)\ dt - \omega^2  \int\limits_{0}^{a}t^{2n-3}\alpha(t)\ dt$$
$$ \le\frac{\omega^2}{a^2} \int\limits_0^{b} t^{2n-1} \alpha(t)\ dt -  \omega^2 \int\limits_{0}^{b}t^{2n-3}\alpha(t)\ dt,
$$
provided all the integrals exist.
\end{lemma}
\pf By a simple rearrangement of integrals, the inequality follows from
$$a^2\int_a^b t^{2n-3} \alpha(t)\ dt \le \int_a^b t^{2n-1} \alpha(t)\ dt. \qed$$

Denote by
$$d_{G}(K,L)=\inf \{b/a: a,b >0 \mbox{  and } \ aK \subset L \subset bK \}$$ the geometric distance between 
two origin-symmetric convex bodies $L$ and $K$ in $\R^{2n}.$
 For a complex star body $K$ in $\R^{2n}$ denote by
$$
d_{IC}(K)=\min \{d_{G}(K, L): \mbox{ $L$ is a complex intersection body in $\R^{2n}$}\}.
$$

\begin{theorem}\label{isom-complex} Let $K$ and $M$ be origin symmetric complex star bodies in $\R^{2n},$
and let $\mu$ be a measure on $\R^{2n}$with even continuous non-negative density $f.$ 
Suppose that for every $\xi\in \SP^{2n-1}$
\begin{equation}\label{stab-ineq}
\mu(K\cap H_\xi) \le \mu(M\cap H_\xi).
\end{equation}
Then 
$$\mu(K)\le (d_{IC}(K))^2\ \mu(M).$$
\end{theorem}

\pf  Without loss of generality, we can assume that the density $f$ is invariant with respect to rotations
$R_\theta.$ In fact, we can consider the measure $\mu_c$ with the density 
$$f_c(x)=\frac1{2\pi}\int_0^{2\pi} f(R_\theta(x))d\theta,$$ 
then $\mu_c(K\cap H_\xi)=\mu(K\cap H_\xi)$ and $\mu_c(K)=\mu(K)$ 
for any complex star body $K$ in $\R^{2n}$ and any $\xi\in \SP^{2n-1}.$

By (\ref{meas-sect-polar}), the condition (\ref{stab-ineq}) can be written as 
$${\cal{R}}_c \left(\int_0^{\|\cdot\|_K^{-1}} r^{2n-3} f(r\cdot)\ dr \right)(\xi) $$
\begin{equation}\label{eq-radon}
\le  {\cal{R}}_c \left(\int_0^{\|\cdot\|_L^{-1}} r^{2n-3} f(r\cdot)\ dr \right)(\xi), \qquad \forall \xi\in \SP^{2n-1}.
\end{equation}
Let $L$ be a complex intersection body in $\R^{2n}$ such that $L\subset K \subset  d_{IC}(K) L.$ Integrate (\ref{eq-radon}) over $\SP^{2n-1}$ with respect to the
measure $\mu$ corresponding to the intersection body $L$ by (\ref{defcompint}). By (\ref{defcompint}) 
$$
\int\limits_{\SP^{2n-1}}\|x\|_L^{-2}\int\limits_{0}^{\|x\|^{-1}_K} t^{2n-3}
 f(tx)dt\ dx $$
\begin{equation}\label{eq11}
\le \!\!\int\limits_{\SP^{2n-1}}\|x\|_L^{-2} \int\limits_{0}^{\|x\|^{-1}_M} t^{2n-3} f (tx) dt\ dx.
\end{equation}
By Lemma \ref{elemcomp} with $\omega =\|x\|_L^{-1}, a =\|x\|_K^{-1}, b=\|x\|_M^{-1}$ and $\alpha(t)=f(tx),$
$$
\frac{\|x\|_L^{-2}}{\|x\|_K^{-2}} \int\limits_0^{\|x\|_K^{-1}} t^{2n-1} f(tx)dt -\|x\|_L^{-2}  \int\limits_{0}^{\|x\|_K^{-1}}t^{2n-3} f(tx)dt$$
\begin{equation}\label{eq12} 
\le \frac{\|x\|_L^{-2}}{\|x\|_K^{-2}} \int\limits_0^{\|x\|_M^{-1}} t^{2n-1} f(tx)dt -\|x\|_L^{-2}  \int\limits_{0}^{\|x\|_M^{-1}}t^{2n-3} f(tx)dt.
\end{equation}
Integrating (\ref{eq12}) over $\SP^{2n-1}$ and adding it to (\ref{eq11}) we get  
\begin{equation}
\int_{\SP^{2n-1}}\frac{\|x\|_L^{-2}}{\|x\|_K^{-2}} \int\limits_0^{\|x\|_K^{-1}} t^{2n-1} f(tx)dt\ dx 
 \le \int_{\SP^{2n-1}}\frac{\|x\|_L^{-2}}{\|x\|_K^{-2}} \int\limits_0^{\|x\|_M^{-1}} t^{2n-1} f(tx)dt\ dx.
\end{equation}
Since $L\subset K \subset  d_{IC}(K) L,$ the latter inequality gives
$$
\frac{1}{(d_{IC}(K))^2}\int_{\SP^{2n-1}} \int\limits_0^{\|x\|_K^{-1}} t^{2n-1} f(tx)dt\ dx 
 \le \int_{\SP^{2n-1}}\int\limits_0^{\|x\|_M^{-1}} t^{2n-1} f(tx)dt\ dx.
$$
The result follows from (\ref{meas-polar}). \endpf

\begin{corollary}  \label{main2} Suppose that $K$ and $M$ are origin-symmetric complex convex bodies in $\R^{2n}$ 
and $\mu$ is an arbitrary measure on $\R^{2n}$ with even continuous density so that 
$$\mu(K\cap H_\xi) \le \mu(M\cap H_\xi),\qquad \forall \xi\in \SP^{2n-1},$$
then
$$\mu(K)\le 2n\ \mu(M).$$
\end{corollary}

\pf By John's theorem  (see, for example, \cite[Section 3]{MS} or \cite[Theorem 4.2.12]{Ga}), there exists an origin symmetric ellipsoid ${\cal{E}}$ such that
$$\frac 1{\sqrt{2n}} {\cal{E}} \subset K \subset {\cal{E}}.$$ 
Construct a new body ${\cal{E}}_c$ by
$$\|x\|_{{\cal{E}}_c}^{-2}= \frac 1{2\pi} \int_0^{2\pi} \|R_\theta x\|_{\cal{E}}^{-2} d\theta.$$
Clearly, ${\cal{E}}_c$ is $R_\theta$-invariant, so it is a complex star body. For every $\theta\in [0,2\pi]$ the 
distribution $\|R_\theta x\|_{\cal{E}}^{-2}$ is positive definite, because of the connection between the Fourier
transform and linear transformations. So $\|x\|_{{\cal{E}}_c}^{-2}$ is
also a positive definite distribution, and, by \cite[Theorem 4.1]{KPZ}, 
${\cal{E}}_c$ is a complex intersection body. Since $\frac 1{\sqrt{2n}} {\cal{E}}\subset K \subset {\cal{E}}$
and $K$ is $R_\theta$-invariant as a complex convex body, we have 
$$\frac 1{\sqrt{2n}} R_\theta {\cal{E}} \subset K \subset R_\theta {\cal{E}}, \quad \forall \theta\in [0,2\pi],$$
so
$$\frac 1{\sqrt{2n}} {\cal{E}}_c\subset K \subset {\cal{E}}_c.$$
Therefore, $d_{IC}(K)\le \sqrt{2n},$ and the result follows from Theorem \ref{isom-complex}. \endpf

\section{The case of convex measures}\label{logconcave}

Following works of  Borell \cite{Bor1, Bor2}, we define the classes of $s$-concave measures.  Let $s \in [-\infty, 1]$. A measure $\mu$ on $\R^n$ is called {\it $s$-concave} if for any compact $A,B \subset \R^n$, with $\mu(A)\mu(B) >0$
and $0 < \lambda < 1$, we have
$$
\mu(\lambda A + (1-\lambda)B) \ge  (\lambda\mu(A)^s+(1-\lambda)\mu(B)^{s})^{1/s}.
$$
 The case where $s=0$ corresponds to {\it log-concave} measures
$$
\mu(\lambda A + (1-\lambda)B) \ge  \mu(A)^\lambda \mu(B)^{(1-\lambda)},
$$
and the case $s=-\infty$ corresponds to {\it convex measures}:
$$
\mu(\lambda A + (1-\lambda)B) \ge  \min\{\mu(A),\mu(B)\}.
$$
We also note that the class of convex measures is the largest class in this group in the sense that it contains 
all other $s$-concave measures. Due to this fact, we concentrate our attention on convex measures.

Borell \cite{Bor1, Bor2} has shown that a measure $\mu$ on $\R^n$ whose support is not contained
in any affine hyperplane is a convex measure if and only if it is
absolutely continuous with respect to Lebesgue measure, and its
density $f$ is a {\it $-1/n$-concave function} on its support, i.e.
$$
f(\lambda x + (1-\lambda) y) \ge (\lambda f(x)^{-1/n}+(1-\lambda) f(y)^{-1/n})^{-n}
$$
for all $x,y: f(x), f(y) >0$ and $\lambda \in [0,1]$. Note  that it follows from the latter definition that if $f(x)$ is a $-1/n$-concave function then $f^{-1/n}$ is a convex function on its support.

We need the following theorem of Bobkov (\cite{Bob}, Theorem 2.1) which is a generalization of the previous result of  Ball \cite{Ba2} (we also refer to \cite{CFPP} for a simpler proof).

\begin{theorem} \label{th:ball}
Let $f:\R^n\rightarrow [0,\infty)$ be an even $-1/n$-concave function on its support, satisfying $0<\int_{\R^n} f <\infty$. Then the map
$$
x \rightarrow  \left(\int_0^\infty f(rx)r^{n-2}dr\right)^{-\frac{1}{n-1}}
$$
defines a norm on $\R^n$.
\end{theorem}

An immediate consequence of the Ball-Bobkov theorem is a  very useful technique of connecting a convex measure 
of one convex body with volume of another convex body. This techniques allows to generalize a number of classical results on Lebesque measure to the case of convex measures  (see \cite{Ba2}, \cite{Bob}, \cite{KYZ} and \cite{CFPP}). Namely, given the density $f$ of a convex measure $\mu$ and a convex symmetric body $K$ we define a body $K_f$ by
$$
\|x\|_{K_f} = \left((n-1)\int_0^\infty (1_K f)(rx)r^{n-2}dr\right)^{-\frac{1}{n-1}},
$$ 
where $1_K$ is the indicator function of $K.$ Theorem  \ref{th:ball} guarantees that $K_f$ is convex. Moreover, by (\ref{volume=spherradon})
\begin{align}
\vol_{n-1}(K_f \cap \xi^\perp)=&\frac{1}{n-1}\int\limits_{\SP^{n-1} \cap \xi^\perp}\|x\|_{K_f} ^{-(n-1)}dx\nonumber \\ =&\int\limits_{\SP^{n-1} \cap \xi^\perp}\int_0^\infty (1_K f)(rx)r^{n-2}dr \ dx=\mu(K\cap \xi^\perp).
 \end{align}
 Our next goal is to estimate $\vol_n(K_f)$. We start with a lemma on the behavior of $-1/n$-concave functions,  the proof of which may be found in  \cite[Lemma 2.4]{Kl}  and \cite[Lemma 4.2]{Bob}.
 \begin{lemma}\label{lm:klartag} Let $n\ge 1$ be an integer, and let $g:[0,\infty) \to [0,\infty)$ be a $-1/n$-concave, non-increasing 
function with $g(0) = 1$, $0 <\int_{0}^\infty  g(t) t^{n-2} dt < \infty$. Then
 $$
 c_1 \le \frac{\int_0^\infty t^{n-1} g(t) dt}{\left(\int_0^\infty t^{n-2} g(t) dt\right)^{\frac{n}{n-1}}} \le c_2, $$
 where $c_1, c_2 >0$ are universal constants.
 \end{lemma}
 \noindent{\bf Remark:} We need the $-1/n$-concavity assumption only to prove the right-hand side inequality in the above lemma. The left-hand side does not require this assumption, but does require $g\le e^n$ .
\smallbreak

Now assume that $f(0)=1,$ $f$ is even and $-1/n$-concave, then $f(tx)$ is non-increasing for $t \ge 0$ and the
function  $g(t)= (1_K f)(tx)$ satisfies the conditions of Lemma \ref{lm:klartag}. 
By (\ref{polar-volume})
$$
\vol_n(K_f)=\frac{1}{n}\int_{\SP^{n-1}} \|x\|_{K_f}^{-n} dx =\frac{1}{n}\int_{\SP^{n-1}}    \left((n-1)\int_0^\infty (1_K f)(rx)r^{n-2}dr\right)^{\frac{n}{n-1}} dx
$$
and applying Lemma \ref{lm:klartag} we get
\begin{equation}
c_1\mu(K)\le\vol_n(K_f) \le c_2\mu(K),
\end{equation}
where $c_1, c_2 >0$ are universal constants (and the right-hand side inequality does not require $-1/n$-concavity, but does require boundness of $f$).

We refer to \cite{MP} for the definition of the isotropic constant $L_K$ of a convex body $K$. It was proved in \cite{MP} that if $L_n=\max\{L_K: \mbox{K is a convex symmetric body  in }\ \R^n\}$ 
then from
$$
\vol_{n-1}(K\cap \xi^\perp) \le  \vol_{n-1}(M\cap \xi^\perp)
$$
we get  $\vol_n(K) \le cL_n \vol_n(M)$. Applying this fact to bodies $K_f$ and $M_f$ we immediately get the following theorem.
\begin{theorem}\label{th:log}  For  any measure $\mu$ with continuous, non-negative even $-1/n$-concave density $f$ on $\R^n$ and any  two
convex origin-symmetric bodies $K, M \subset \R^n$ such
that
\begin{equation}\label{e:bnlog}
\mu(K\cap \xi^\bot)\le\mu(M\cap \xi^\bot),\ \qquad \forall \xi\in \SP^{n-1}
\end{equation}
we have
$$
\mu(K)\le   cL_n \mu(M).
$$
\end{theorem}
\noindent{\bf Remark:} We note that the assumption $f(0)=1$ is not necessary in the above theorem due to the fact that the theorem does not change when $\mu$ is multiplied by a constant.
It was proved by Bourgain that $L_n \le c n^{1/4} \log(n+1)$ and the $\log(n+1)$ factor was after removed by Klartag \cite{Kl}, which implies the following corollary.

\begin{corollary} For  any convex measure $\mu$ with continuous, non-negative even density $f$ on $\R^n$ and any  two
convex origin-symmetric bodies $K, M \subset \R^n$ such
that
\begin{equation}
\mu(K\cap \xi^\bot)\le\mu(M\cap \xi^\bot),\ \qquad \forall \xi\in \SP^{n-1}
\end{equation}
we have
$$
\mu(K)\le   c n^{1/4} \mu(M).
$$
\end{corollary} 

It was also proved in \cite{MP} that for any convex origin-symmetric body $K\subset \R^n$
$$
\max_{\xi \in \SP^{n-1}} \vol_{n-1}(K \cap \xi^\perp) \ge \frac{c}{L_K} \vol_n(K)^{\frac{n-1}{n}}.
$$
which gives (applying the latter inequality to $K_f$)
$$
\frac{1}{f(0)}\max_{\xi \in \SP^{n-1}} \mu(K \cap \xi^\perp) \ge \frac{c}{L_K} \left(\frac{\mu(K)}{f(0)}\right)^{\frac{n-1}{n}},
$$
which implies
\begin{corollary}\label{co:bob} For  any convex measure $\mu$ with continuous, non-negative even density $f$ on $\R^n$ and any  
convex origin-symmetric  body $K \subset \R^n$ we have
$$
\max_{\xi \in \SP^{n-1}} \mu(K \cap \xi^\perp) \ge \frac{c}{L_K} \mu(K)^{\frac{n-1}{n}} f(0)^{\frac{1}{n}}.
$$
\end{corollary}
Using convexity of $\mu$ we get that $\frac{\mu(K)}{f(0)} \le \vol_n(K)$, which proves the following hyperplane inequality 
for convex measures.
\begin{corollary} \label{co: koldobsky} For  any convex measure $\mu$ with continuous, non-negative even density $f$ on $\R^n$ and any  
convex origin-symmetric  body $K \subset \R^n$ we have
$$
\max_{\xi \in \SP^{n-1}} \mu(K \cap \xi^\perp) \ge \frac{c}{L_K} \mu(K) \vol_n(K)^{-\frac{1}{n}},
$$
and thus
$$
\mu(K)\le Cn^{1/4} \max_{\xi \in \SP^{n-1}} \mu(K \cap \xi^\perp)\ \vol_n(K)^{\frac{1}{n}}.
$$
\end{corollary}

We would like to note that Corollary \ref{co:bob} was essentially proved by Bobkov  \cite[Theorem 4.1]{Bob}.
Our goal is a generalization
of the hyperplane inequality to the case of most general measures with positive even and continuous density.  We note that Corollary \ref{co:bob} is false in the case of general measures. Indeed, consider $f(x)=1/(1+|x|^p)$ for some $p\in (0, n)$, then $f(x)$ is radial decreasing and  $f(0)=1$ is still the maximum for $f$ on $\R^n$. Let $K=tB_2^n$ for $t$ large enough, then using (\ref{polar-measure}) we get
$$
\mu(tB_2^n)=\frac{1}{|\SP^{n-1}|}\int_0^t \frac{r^{n-1}dr}{1+r^p} \ge \frac{c t^{n-p}}{(n-p)|\SP^{n-1}|}
$$
and
$$
\mu(tB_2^{n-1})\le \frac{ t^{n-p-1}}{(n-p-1)|\SP^{n-2}|}.
$$
Thus for Corollary \ref{co:bob}  to be correct we must have for all large $t$
$$
t^{n-p-1} \ge c_n\  t^{(n-p)\frac{n-1}{n}}
$$
or
$$
n-p-1 \ge (n-p)\frac{n-1}{n} \qquad {\rm and} \qquad  p\le 0,
$$
which gives a contradiction. 

We finish this note with an observation related to the hyperplane inequality for measures.
\begin{lemma}\label{lm:const} For  any measure $\mu$ with continuous, non-negative even density $f$ on $\R^n$ consider a symmetric star-shaped body $K \subset \R^n$ such that
$\mu(K \cap \xi^\perp) =\mu(K \cap \theta^\perp)$ for all $\xi, \theta \in \SP^{n-1}$, then 
$$
 \mu(K) \le C \mu(K \cap \xi^\perp)\ \vol_n(K)^{\frac{1}{n}},\qquad \forall \xi\in \SP^{n-1}.
 $$
\end{lemma}

\beginproof Assume $\mu(K \cap \xi^\perp)= \Lambda$, then applying (\ref{measure=spherradon}) we get
\begin{equation}\label{eq:equalsec}
\RRR\left(\int_0^{\|\cdot\|_K^{-1}} r^{n-2}f(r\ \cdot)\ dr \right)(\xi) = \Lambda,
\end{equation}
and applying the Funk-Minkowski uniqueness theorem for the spherical Radon transform (see for example \cite{Kbook}) we get
$$
\int_0^{\|\xi\|_K^{-1}} r^{n-2}f(r \xi )\ dr = \frac{ \Lambda}{ \vol_{n-2}(\SP^{n-2})}, \qquad  \forall \xi \in \SP^{n-1}.
$$
We also note that
$$
\int_0^{\|\xi\|_K^{-1}} r^{n-1}f(r \xi)\ dr \le \|\xi\|_K^{-1} \int_0^{\|\xi\|_K^{-1}} r^{n-2}f(r \xi)\ dr = \frac{ \Lambda \|\xi\|_K^{-1}}{ \vol_{n-2}(\SP^{n-2})}.
$$
Finally, integrating the above inequality over $\xi \in \SP^{n-1}$ and  applying (\ref{polar-measure})
\begin{align}
\mu(K) \le & \frac{ \Lambda}{ \vol_{n-2}(\SP^{n-2})} \int_{\SP^{n-1}}  \|\xi\|_K^{-1} d \xi \\ \nonumber  \le  & \frac{ \Lambda \vol_{n-1}(\SP^{n-1})^{\frac{n-1}{n}}}{ \vol_{n-2}(\SP^{n-2})}
\left(\int_{\SP^{n-1}}  \|\xi\|_K^{-n} d \xi \right)^{1/n} \\ \nonumber \le & C \Lambda \vol_n(K)^{1/n}.
\end{align}
\endproof
\noindent{\bf Remark:}  We note that the body $K$ in Lemma \ref{lm:const} exists for all $\Lambda>0$ such that
$$
\Lambda \le  \vol_{n-2}(\SP^{n-2})  \min_{\xi \in \SP^{n-1}}  \int_0^{\infty} r^{n-2}f(r \xi)\ dr.
$$
This follows from (\ref{eq:equalsec}), properties of the spherical Radon transform and the fact that $f(x) \ge 0$ (see Corollary 1 in \cite{Z2}). Clearly, $K$ is not necessarily a convex body. It seems to be quite difficult to find a sufficient condition on $f$ for $K$ to be convex. For any rotation invariant $f$ we get that $K$ is a dilate of the Euclidean ball.

 \end{document}